\documentclass{amsart}
\theoremstyle{definition}

\theoremstyle{remark}

\numberwithin{equation}{section}

\begin{document}

\title{Nonnegative rank of a matrix\\ with one negative eigenvalue}

\author{Yaroslav Shitov}
\address{National Research University Higher School of Economics, 20 Myasnitskaya Ulitsa, Moscow 101000, Russia}
\email{yaroslav-shitov@yandex.ru}


\begin{abstract}
We show that a rank-three symmetric matrix with exactly one negative eigenvalue can have arbitrarily large nonnegative rank.
\end{abstract}

\maketitle

The concept of \textit{nonnegative rank}, equal to the smallest integer $r$
for which a given nonnegative matrix can be written as a product
of an $m$-by-$r$ and an $r$-by-$n$ nonnegative matrices, provides a useful
tool for data analysis, statistics, complexity theory and other applications.
Problems related to upper and lower bounds for the nonnegative rank take an
important place, and we solve one of those problems posed by Beasley and Laffey
in~\cite{BL} by proving the result mentioned in the abstract.

Fix an integer $k\geq3$ and take $h_{i+1}=\mathrm{e}^{-1/h_i}$,
for $i\leq k$ and a small real $h_1>0$.
Define the $3$-by-$3$ matrix $\Lambda$ by taking
$\Lambda_{rs}=-1$ if $\{r,s\}=\{2,3\}$ and $\Lambda_{rs}=1$ otherwise;
the $2k$-by-$3$ matrix $B$ by $B_{2i-1,1}=B_{2i,1}=1$,
$B_{2i,2}=14 + h_i$, $B_{2i-1,3}=6 + h_i$,
$B_{2i-1,2}=4 + h_i - \frac{h_i}{7 + h_i^2}$ and
$B_{2i,3}=\frac{21 + 8 h_i + 4 h_i^2 - h_i^3}{21 + h_i + 3 h_i^2}$.
Then the matrix $A=B\Lambda B^\top$ is symmetric and has rank $3$,
and every principal $3$-by-$3$ minor of $A$ is negative as $\det\Lambda<0$.
So the characteristic polynomial $\det(A-\lambda I)$ has either
one or three negative roots. The latter is impossible due to the Perron--Frobenius
theorem as we have $A_{pq}=0$ if $\{p,q\}=\{2i-1,2i\}$ for some $i$ and $A_{pq}>0$ otherwise.
The nonnegative rank of $A$ exceeds $\log_2 k$ by Proposition~1.3 from~\cite{BL}.

\end{document}